\documentclass[a4paper,11pt,oneside,openany]{article}

\usepackage{amsmath,amssymb}
\usepackage{amsfonts}

\usepackage{graphics}
\usepackage{float}
\usepackage[dvips]{graphicx,epsfig}
\usepackage{subfigure}
\usepackage{latexsym,euscript,epic,eepic}
\newcommand{\proba}{{\mathbb{P}}}
\newcommand{\R}{{\mathbb{R}}}

\newcommand{\Ind}{{\mathbb{I}}}

\newcommand{\WH}{W^{\mbox{\rm\tiny H}}}
\newcommand{\WZ}{W^{\mbox{\rm\tiny Z}}}

\newcommand{\egloi}{\stackrel{d}{=}}
\newcommand{\egdef}{\stackrel{def}{=}}
\newcommand{\eqproba}{\stackrel{{\rm{P}}}{\sim}}
\newcommand{\cvproba}{\stackrel{{\rm{P}}}{\rightarrow}}
\newcommand{\cvloi}{\stackrel{d}{\rightarrow}}

\newcommand{\CQFD}
{%
\mbox{}%
\nolinebreak%
\hfill%
\rule{2mm}{2mm}%
\medbreak%
\par%
}

\newtheorem{Theo}{Theorem}
\newtheorem{Prop}{Proposition}
\newtheorem{Coro}{Corollary}
\newtheorem{Lem}{Lemma}
\DeclareMathOperator*{\argmin}{arg\,min}

\begin{document}

\begin{center} {\bf{FUNCTIONAL NONPARAMETRIC ESTIMATION OF CONDITIONAL EXTREME QUANTILES }}

\vspace*{6mm}
Laurent GARDES\footnote{Corresponding author}, St\'ephane GIRARD and Alexandre LEKINA\\
\vspace*{6mm}
{\small \noindent Team Mistis, INRIA Rh\^{o}ne-Alpes and LJK,\\
655, avenue de l'Europe, Montbonnot\\
38334 Saint-Ismier cedex, France. \\
{\tt Laurent.Gardes@inrialpes.fr}}
\end{center}

\vspace*{3mm}
\ \\
\noindent {\bf{Abstract}} $-$ We address the estimation of quantiles from heavy-tailed distributions when functional covariate information is available 
and in the case where the
order of the quantile converges to one as the sample size increases.
Such "extreme" quantiles
can be located in the range of the data or near and even beyond the boundary of the sample, depending on the convergence rate of their order to one. 
Nonparametric estimators of these functional extreme quantiles are introduced,
their asymptotic distributions are established and their finite sample behavior is investigated. \\ 


\vspace*{5mm}
\ \\
\section{Introduction}

An important literature is dedicated to the estimation of extreme quantiles,
{\it i.e.} quantiles of order $1-\alpha$ with $\alpha$ tending to zero.
The most popular estimator was proposed by Weissman~\cite{wei78},
in the context of heavy-tailed distributions, and
adapted to Weibull-tail distributions in~\cite{DGGG2,CSTM}. 
We also refer to~\cite{Dek} for the general case. \\
In a lot of applications, some covariate information is recorded simultaneously with the quantity of interest. For instance, in climatology one may be interested in the estimation of return periods associated to extreme rainfall as a function of the geographical location. 
The extreme quantile thus depends on the covariate and is referred in the sequel to as the conditional extreme quantile.  
Parametric models for conditional extremes are proposed
in~\cite{davsmi90,smith89} whereas semi-parametric methods are
considered in~\cite{beigoe03,haltaj00}.
Fully non-parametric estimators have been first introduced in~\cite{davram00}, 
where a local polynomial modelling
of the extreme observations is used. Similarly, spline estimators
are fitted in~\cite{chadav05} through a penalized maximum likelihood method.
In both cases, the authors focus on univariate covariates and on the
finite sample properties of the estimators.
These results are extended in~\cite{beigoe04} where local polynomials estimators
are proposed for multivariate covariates and where their asymptotic properties
are established.

Besides, covariates may be curves in many situations coming from
applied sciences such as chemometrics (see Section~\ref{realdata} for an
illustration) or astrophysics~\cite{StaCom}.
However, the estimation of conditional extreme quantiles
with functional covariates has not been addressed yet.
Two statistical fields are involved in this study.
In the one hand, nonparametric smoothing techniques adapted to functional data
are required in order to deal with the covariate. 
We refer to~\cite{Bosq,fer06,Ramsay,Ramsay2}
for overviews on this literature.
We propose here to select the observations to be used in the conditional quantile estimator by a moving window approach. In the second hand, once this selection is achieved, extreme-value methods are used to estimate the conditional quantile, see~\cite{EMBR} for a comprehensive treatment
of extreme-value methodology in various frameworks.

Whereas no parametric assumption is made on the functional covariate,
we assume that the conditional distribution is heavy-tailed.
This semi-parametric assumption amounts to supposing that the conditional survival function decreases at a polynomial rate.
To estimate the conditional quantile, we focus on three different situations. In the first one, the convergence of $\alpha$ to zero is slow enough so that the quantile is located in the range of the data. In the second situation,
the quantile is located near the boundary of the sample.
Finally, in the third situation, the convergence of $\alpha$ to zero is sufficiently fast so that the quantile may be beyond the boundary of the sample. This situation is clearly the most difficult one since an extrapolation outside the range of the sample is needed to achieve the estimation. 

Nonparametric estimators are defined in Section~\ref{defest} for each situation.
Their asymptotic distributions are derived in Section~\ref{secmain}.
Some examples are provided in Section~\ref{examples} and an illustration on spectrometric data is given in Section~\ref{realdata}. Proofs are postponed to Section~\ref{secpreuves}.

\section{Estimators of conditional extreme quantiles}
\label{defest}

Let $E$ be a (finite or infinite dimensional) metric space associated to a metric $d$. Let us denote by $F(.,x)$ the conditional cumulative distribution function of a real random variable $Y$ given $x \in E$ and by $q(\alpha,x)$ the associated conditional quantile of order $1-\alpha$ defined by 
$$
F(q(\alpha,x),x)=1-\alpha,
$$
for all $x\in E$ and $\alpha\in(0,1)$.
In this paper, we focus on the case where, for all $x \in E$, $F(.,x)$ is the cumulative distribution function of a heavy-tailed distribution. 
In such a situation, the conditional quantile $q(.,x)$ satisfies, for all $\lambda>0$, 
\begin{equation}
\label{modelq}
\lim_{\alpha \to 0} \frac{q(\lambda\alpha,x)}{q(\alpha,x)} = \lambda^{-\gamma(x)},
\end{equation}
where $\gamma(.)$ is an unknown positive function of the covariate $x$
referred to as the conditional tail index.
Loosely speaking, the conditional quantile $q(.,x)$ decreases towards 0 at 
a polynomial rate driven by $\gamma(x)$.
The conditional quantile 
is said to be regularly varying at 0 with index $-\gamma(x)$,
and this property characterizes heavy-tailed distributions.
We refer to~\cite{BING} for a general account on regular variation theory
and to paragraph~\ref{exlois} for some examples of distributions
satisfying~(\ref{modelq}).

Given a sample $(Y_1,x_1), \ldots, (Y_n,x_n)$ of independent observations, our aim is to build point-wise estimators 
of conditional quantiles. More precisely, for a given $t\in E$, we want to estimate $q(\alpha,t)$, focusing on the case where the design points $x_1,\ldots,x_n$ are non random.
To this end, for all $r>0$, let us denote by $B(t,r)$
the ball centered at point $t$ and with radius $r$ defined by
$$
B(t,r)=\{ x \in E, \ d(x,t) \leq r \}
$$
and let $h_{n,t}=h_t$ be a positive sequence tending to zero as $n$ goes to infinity. The proposed estimator uses a moving window approach since it is based on the response variables $Y_i's$ for which the associated covariates $x_i's$ belong to the ball $B(t,h_t)$.
The proportion of such design points is thus defined by
$$
\varphi(h_t)=\frac{1}{n}\sum_{i=1}^n \Ind
\{ x_i\in B(t,h_t) \}
$$
and plays an important role in this study.
It describes how the design points concentrate in the neighborhood of
$t$ when $h_t$ goes to zero, similarly to the small ball probability does,
see for instance the monograph on functional data analysis~\cite{fer06}.
Thus, the nonrandom number of observations in the slice $S_t=(0,\infty) \times B(t,h_t)$ is given by $m_{n,t} = m_t = n \varphi(h_t)$. Let $\{ Z_i(t), \ i=1,\ldots,m_t\}$ be the response variables $Y_i's$ for which the associated covariates $x_i's$ belong to the ball $B(t,h_t)$ and let $Z_{1,m_t}(t) \leq \ldots \leq Z_{m_t,m_t}(t)$ be the corresponding order statistics. 

In this paper, we focus on the estimation  of conditional "extreme" quantile of order $1-\alpha_{m_t}$. Here, the word "extreme" means that $\alpha_{m_t}$ tends to zero as $n$ goes to infinity, making kernel based estimators~\cite{fer06bis} non adapted.
In the sequel, three situations are considered:
\begin{itemize}
\item[{\bf{(S.1)}}] $\alpha_{m_t} \to 0$ and $m_t \alpha_{m_t} \to \infty$,
\item[{\bf{(S.2)}}] $\alpha_{m_t} \to 0$, $m_t \alpha_{m_t} \to c \in [1,\infty)$ and $\lfloor m_t \alpha_{m_t} \rfloor \to \lfloor c \rfloor $.
\item[{\bf{(S.3)}}] $\alpha_{m_t} \to 0$ and $m_t \alpha_{m_t} \to c\in [0,1)$,
\end{itemize}
where $\lfloor x \rfloor$ denotes the largest integer smaller than $x$.
Let us highlight that, in the unconditional case,
situations {\bf{(S.1)}} and {\bf{(S.3)}} with $c\neq 0$ have already been examined by Dekkers and de Haan~\cite{Dek}, 
the extreme case $c=0$ being considered in~\cite{deH}, Theorem~5.1.
A summary of their results can be found in~\cite{EMBR}, Theorem~6.4.14 and
Theorem~6.4.15.
In situation {\bf{(S.1)}}, $\alpha_{m_t}$ goes to 0 slower than $1/m_t$ and the point-wise estimation of the conditional extreme quantile relies on an interpolation inside the sample, since,
from Proposition~\ref{compamax} below, $q(\alpha_{m_t},t)$ is
eventually almost surely
smaller that the maximal observation $Z_{m_t,m_t}(t)$ in the slice $S_t$. In such a situation, we propose to estimate $q(\alpha_{m_t},t)$ by:
\begin{equation}
\label{q1}
{\hat{q}}_1(\alpha_{m_t},t) = Z_{m_t-\lfloor m_t\alpha_{m_t} \rfloor +1,m_t}(t). \end{equation}
In the intermediate situation {\bf{(S.2)}},
estimator~(\ref{q1}) can still be used, since for $n$ large enough,
$\lfloor m_t\alpha_{m_t} \rfloor = \lfloor c\rfloor >0$ and
thus the estimation relies on a conditional extreme value of the sample.
Let us note that, if $c$ is not an integer, then $m_t \alpha_{m_t} \to c$ implies $\lfloor m_t \alpha_{m_t} \rfloor \to \lfloor c \rfloor$. Otherwise,
if $c$ is an integer, then condition $\lfloor m_t \alpha_{m_t} \rfloor \to \lfloor c \rfloor$ is necessary to prevent the sequence $\lfloor m_t \alpha_{m_t} \rfloor$ from having two adherence values and ${\hat{q}}_1(\alpha_{m_t},t)$
from oscillating.
In situation {\bf{(S.3)}}, $\alpha_{m_t}$ goes to 0 at
the same speed or faster than $1/m_t$ and the conditional extreme quantile is eventually larger than $Z_{m_t,m_t}(t)$ with positive probability 
${\rm{e}}^{-c} \geq {\rm{e}}^{-1}$. Thus, its estimation is more difficult since it
requires an estimation outside the sample. We propose in this case to estimate $q(\alpha_{m_t},t)$ by:
\begin{eqnarray}
\label{q2}
{\hat{q}}_2(\alpha_{m_t},t) & = & {\hat{q}}_1(\beta_{m_t},t) \left ( \beta_{m_t} / \alpha_{m_t} \right )^{{\hat{\gamma}}_n(t)} \nonumber \\
\ & = & Z_{m_t-\lfloor m_t\beta_{m_t} \rfloor +1,m_t}(t)\left ( \beta_{m_t} /\alpha_{m_t} \right )^{{\hat{\gamma}}_n(t)},
\end{eqnarray}
where $\beta_{m_t}$ satisfies {\bf{(S.1)}} and ${\hat{\gamma}}_n(t)$ is a point-wise estimator of the conditional tail index $\gamma(t)$. Such estimators have been proposed both in the finite dimensional setting~\cite{beigoe04} and in the general case~\cite{JMVA}, see also paragraph~\ref{exestis} for some examples. Note that~(\ref{q2}) is an adaptation of Weissman estimator~\cite{wei78} in the case where covariate information is available. The extrapolation is achieved
thanks to the multiplicative term 
$\left ( \beta_{m_t} /\alpha_{m_t} \right )^{{\hat{\gamma}}_n(t)}$ which
magnitude is driven by the estimated tail index ${\hat{\gamma}}_n(t)$.
As expected, the extrapolation is all the more important as the tail is heavy.

\section{Main results}
\label{secmain}

\noindent We first give some notations and conditions useful to establish the asymptotic distributions of our estimators. In the sequel, we fix $t \in E$
and we assume:
\begin{itemize}
\item[{\bf{(A)}}] The conditional quantile function
$$
\alpha\in(0,1)\mapsto q(\alpha,t) \in (0,+\infty)
$$
is differentiable, 
the function defined by 
$$
\alpha\in(0,1)\mapsto \Delta(\alpha,t) = \gamma(t)+\alpha \frac{\partial \log q} {\partial \alpha} (\alpha,t)\in (0,+\infty)
$$
is continuous and such that
$
\displaystyle\lim_{\alpha \to 0} \Delta(\alpha,t)=0.
$
\end{itemize}
Assumption {\bf{(A)}} controls the behavior of the log-quantile function
with respect to its first variable. It is a sufficient condition to
obtain the heavy-tail property~(\ref{modelq}),
see for instance~\cite{BING}, Chapter~1. 
For all $a \in (0,1)$, let us introduce
$$
{\bar{\Delta}}(a,t)=\sup_{\alpha \in (0,a)} |\Delta(\alpha,t)|.
$$
The largest oscillation of the log-quantile function with respect to its second variable is defined for all $a\in (0,1/2)$ as
\[ 
\omega_n(a) = \sup \left \{ \left | \log \frac{q(\alpha,x)}{q(\alpha,x')} \right | \ {, } \ \alpha \in ( a , 1-a) \ {, } \ (x,x') \in B(t,h_t)^2 \right \}.
\]
Finally, let $k_t \in \{1,\dots, m_t\}$ and $J_{k_t}=\{ 1,\ldots,k_t \}$.
Our first result establishes a representation in distribution
of the largest random variables of the sample $Z_i(t)$, $i\in \{1,\dots,m_t\}$.
\begin{Prop}
\label{proplaurent}
If 
$k_t / m_t \to 0$ and  $k_t^2\omega_{n}(m_t^{-(1+\delta)})\to 0$
for some $\delta>0$, then, there exists an event ${\cal{A}}_n$ with $\proba({\cal{A}}_n) \to 1$ 
as $n\to\infty$ such that 
\[ \left \{ \left ( \log Z_{m_t-i+1,m_t}, i \in J_{k_t} \right ) | {\cal{A}}_n \right \} \egloi \left \{ \left ( \log q(V_{i,m_t},T_i), i \in J_{k_t} \right ) | {\cal{A}}_n  \right \}, \]
where 
$V_{1,m_t} \leq \ldots \leq V_{m_t,m_t}$ are the order statistics associated to the sample $\{V_1,\ldots,V_{m_t}\}$ of independent uniform variables and $\{T_1,\dots,T_{k_t}\}$ are random variables in the ball $B(t,h_t)$.
\end{Prop}
Note that this result is implicitly used in~\cite{JMVA}, proof of Theorem~1.
We also refer to~\cite{FHR}, Theorem~3.5.2,
for the approximation of the
nearest neighbors distribution using the Hellinger distance and to~\cite{Gango}
for the study of their asymptotic distribution.
Here, condition $k_t^2\omega_{n}(m_t^{-(1+\delta)})\to 0$ shows that, 
the smoother the quantile function is on the slice $S_t$, {\it i.e.} the smaller
its oscillation is, the easier the control of the uppest observations is,
{\it i.e} the larger $k_t$ can be.

The next proposition is dedicated to the study of the position
of the conditional extreme quantile $q(\alpha,t)$ with respect to the largest observation in the slice $S_t$.
\begin{Prop}
\label{compamax}
If
$
\omega_{n}(m_t^{-(1+\delta)}) \to 0
$
for some $\delta>0$, then
\begin{itemize}
\item under {\bf{(S.1)}}, $\proba(Z_{m_t,m_t}<q(\alpha_{m_t},t)) \to 0$,
\item under {\bf{(S.2)}} or {\bf{(S.3)}}, $\proba(Z_{m_t,m_t}<q(\alpha_{m_t},t)) \to {\rm{e}}^{-c}$.
\end{itemize}
\end{Prop}
Let us first focus on situation {\bf{(S.1)}}
where the estimation of the conditional extreme quantile is addressed 
using $\hat q_1(\alpha_{m_t},t)$, an upper order statistic chosen in 
the considered slice.
\begin{Theo}
\label{quantclass}
Let $(\alpha_{m_t})$ be a sequence satisfying {\bf{(S.1)}}.\\
If 
$
(m_t \alpha_{m_t})^2\omega_{n}(m_t^{-(1+\delta)}) \to 0
$ 
for some $\delta>0$ then,
\[ (m_t\alpha_{m_t})^{1/2} \left ( \frac{{\hat{q}}_1(\alpha_{m_t},t)}{q(\alpha_{m_t},t)} -1 \right ) \cvloi {\cal{N}}(0,\gamma^2(t)). \]
\end{Theo}
It appears that the estimator is asymptotically Gaussian, with asymptotic 
variance proportional to $\gamma^2(t)/(m_t\alpha_{m_t})$.
Thus, the heavier is the tail, the larger is $\gamma(t)$, and the larger
is the variance. Besides, the asymptotic variance being inversely proportional
to $\alpha_{m_t}$, the estimation remains more stable when the extreme
quantile is far from the boundary of the sample.
Considering now situation {\bf{(S.2)}}, an asymptotically Gaussian behavior
cannot be expected since, in this case, the estimator is based
on the $\lfloor c \rfloor$th uppest order statistic in the considered slice.
\begin{Theo}
\label{quantinter}
Let $(\alpha_{m_t})$ be a sequence satisfying {\bf{(S.2)}}.\\
If 
$
\omega_{n}(m_t^{-(1+\delta)}) \to 0
$ 
for some $\delta>0$ then,
\[ \left ( \frac{{\hat{q}}_1(\alpha_{m_t},t)}{q(\alpha_{m_t},t)} -1 \right ) \cvloi {\cal{E}}(c,\gamma(t)), \]
where ${\cal{E}}(c,\gamma(t))$ is a non-degenerated distribution.
\end{Theo}
The asymptotic distribution ${\cal{E}}(c,\gamma(t))$ could be explicitly deduced from the proof of the result. It is omitted here for the sake of simplicity.
Situation {\bf(S.3)} is more complex since the asymptotic distribution
of $\hat q_2$ may depend both on the behavior of $\hat q_1$ and $\hat\gamma_n$.
In the next theorem, two cases are investigated.
In situation~(i), the asymptotic distribution of $\hat q_2$
is driven by $\hat q_1$. At the opposite, in situation~(ii), 
$\hat q_2$ inherits its asymptotic distribution from $\hat\gamma_n$. 
\begin{Theo}
\label{quantextreme}
Let $(\beta_{m_t})$ be a sequence satisfying {\bf{(S.1)}} 
and let $(\alpha_{m_t})$ be a sequence eventually smaller than $(\beta_{m_t})$.
Define $\zeta_{m_t}=(m_t\beta_{m_t})^{1/2}\log ( \beta_{m_t}/\alpha_{m_t})$.\\
If
 $(m_t\beta_{m_t})^2\omega_n(m_t^{-(1+\delta)}) \to 0$ for some $\delta>0$
and there exists a positive sequence $\upsilon_n(t)$ and a distribution ${\cal{D}}$ such that
\begin{equation}
\label{CVgamma}
\upsilon_n(t) ({\hat{\gamma}}_n(t)-\gamma(t)) \cvloi {\cal{D}},
\end{equation}
then, two situations arise:
\begin{itemize}
\item[(i)]Under the additional condition
\begin{equation}
\label{C1}
\zeta_{m_t} \max \left \{ \upsilon_n^{-1}(t),{\bar{\Delta}}(\beta_{m_t},t) \right \} \to 0,
\end{equation} 
we have
\begin{equation}
\label{cvloi1}
(m_t\beta_{m_t})^{1/2} \left (\frac{{\hat{q}}_2(\alpha_{m_t},t)}{q(\alpha_{m_t},t)} -1\right ) \cvloi {\cal{N}}(0,\gamma^2(t)).
\end{equation}
\item[(ii)] Otherwise, under the additional condition
\begin{equation}
\label{C2}
 \upsilon_n(t) \max \left \{ \zeta_{m_t}^{-1},{\bar{\Delta}}(\beta_{m_t},t) \right \} \to 0,
\end{equation} 
we have
\begin{equation}
\label{cvloi2}
\frac{\upsilon_n(t)}{\log \left ( \beta_{m_t}/\alpha_{m_t} \right )} \left (\frac{{\hat{q}}_2(\alpha_{m_t},t)}{q(\alpha_{m_t},t)} -1\right ) \cvloi {\cal{D}}.
\end{equation}
\end{itemize}
\end{Theo}
Note that, even though the main interest of this result
is to tackle the case where $(\alpha_{m_t})$ is a sequence satisfying
{\bf{(S.3)}}, it can also be applied in the more general situation
where $\alpha_{m_t}$ is eventually smaller than $\beta_{m_t}$. 
For instance, it appears that, in situation {\bf{(S.2)}}, 
${\hat{q}}_2(\alpha_{m_t},t)$ is a consistent estimator 
of $q(\alpha_{m_t},t)$ in the sense that the ratio converges to one
in probability whereas, in view of Theorem~\ref{quantinter},
${\hat{q}}_1(\alpha_{m_t},t)$ is not consistent.
Some applications of Theorem~\ref{quantextreme}
are provided in the next section. 

\section{Examples}
\label{examples}

In paragraph~\ref{exestis}, the above theorem is illustrated
with a particular family of conditional tail index estimators.
The corresponding assumptions are simplified in paragraph~\ref{exlois}
for some classical heavy-tailed distributions.

\subsection{Some conditional tail-index estimators}
\label{exestis}

In~\cite{JMVA}, a family of conditional tail index estimators is introduced.
They are based on a weighted sum of the log-spacings between
the $k_t$ largest order statistics $Z_{m_t-k_t+1,m_t},\dots,Z_{m_t,m_t}$.
The family is defined by 
\begin{equation}
\label{family}
{\hat{\gamma}}_n(t,W) = \sum_{i=1}^{k_{t}} i \log \left ( \frac{Z_{m_{t}-i+1,m_{t}}(t)}{Z_{m_{t}-i,m_{t}}(t)} \right ) W\left(i/k_{t},t\right)  
\left/ \sum_{i=1}^{k_{t}}W\left(i/k_{t},t\right)\right., 
\end{equation}
where $W(.,t)$ is a weight function defined on $(0,1)$ and integrating to one.
Basing on~(\ref{family}) and considering
$\beta_{m,t}=k_t/m_t$,
the conditional extreme quantile estimator~(\ref{q2}) can be written as
$$
{\hat{q}}_2(\alpha_{m_t},t,W) 
  =  Z_{m_t- k_t +1,m_t}(t)\left (\frac{ k_t }{m_t \alpha_{m_t}} \right )^{{\hat{\gamma}}_n(t,W)}.
$$
From~\cite{JMVA}, Theorem~2, under some conditions on the weight function,
${\hat{\gamma}}_n(t,W)$ is asymptotically Gaussian:
$$
k_t^{1/2}({\hat{\gamma}}_n(t,W)-\gamma(t)) \cvloi {\cal{N}}(0,\gamma^2(t){\cal AV}(t,W)), 
$$
where ${\cal AV}(t,W)=\int_0^1W^2(s,t) ds$. Letting $\upsilon_n(t)=k_t^{1/2}$,
we obtain
$$
\zeta_{m_t} \upsilon_n^{-1}(t) = \log \left ( \frac{k_t}{m_t\alpha_{m_t}} \right ) \to \infty,
$$
in situation {\bf (S.2)} or {\bf (S.3)},
which means that condition~(\ref{C1}) cannot be satisfied. Thus, only situation (ii) of~Theorem~\ref{quantextreme} may arise leading to the following corollary.
\begin{Coro}
\label{coro}
Suppose the assumptions of~\cite{JMVA}, Theorem~2 hold. Let $k_t\to\infty$ such that 
\begin{eqnarray}
\label{cond1}
&k_t^{1/2} \bar\Delta(k_t/m_t,t)\to 0 & \mbox{ and } \\
\label{cond2}
&k_t^2     \omega_n(m_t^{-(1+\delta)}) \to 0 &\mbox{ for some } \delta>0.
\end{eqnarray}
Let $(\alpha_{m_t})$ be a sequence satisfying {\bf (S.2)} or {\bf (S.3)}.
Then,
$$
\frac{k_t^{1/2}}{\log ( k_t/(m_t\alpha_{m_t}))} \left (\frac{{\hat{q}}_2(\alpha_{m_t},t,W)}{q(\alpha_{m_t},t)} -1\right ) \cvloi {\cal N}(0,\gamma^2(t){\cal AV}(t,W)).
$$
\end{Coro}
As an example, one can use constant weights $\WH(s,t) = 1$ to obtain the
so-called conditional Hill estimator with ${\cal AV}(t,\WH)=1$ or logarithmic weights $\WZ(s,t)=-\log(s)$
leading to the conditional Zipf estimator with ${\cal AV}(t,\WZ)=2$. We refer to~\cite{JMVA}, Section~4, for further details.

\subsection{Illustration on some heavy-tailed distributions}
\label{exlois}

Standard Pareto distribution is the simplest example of heavy-tailed distribution. Its conditional quantile of order $1-\alpha$ decreases as a power function of $\alpha$ since, in this case,
$
q(\alpha,t)=\alpha^{-\gamma(t)}.
$
Therefore $\Delta(\alpha,t)=0$ for all $\alpha \in (0,1)$ and condition~(\ref{cond1}) of Corollary~\ref{coro} vanishes.
Another example is Fr\'echet distribution for which
$$
q(\alpha,t) = \alpha^{-\gamma(t)} \left\{\displaystyle\frac{1}{\alpha}\log\left(\displaystyle
\frac{1}{1-\alpha}\right)\right\}^{-\gamma(t)}.
$$
Here, the conditional quantile approximatively
decreases as a power function of $\alpha$ since, in this case,
 $q(\alpha,t)\sim\alpha^{-\gamma(t)}$, 
the quality of this approximation being controlled by 
$$
\Delta(\alpha,t) = 
 -\displaystyle\frac{\gamma(t)}{2} \alpha (1+O(\alpha)) \mbox{ as } \alpha\to 0.
$$
A similar example is given by Burr distributions for which
$$
q(\alpha,t) = \alpha^{-\gamma(t)}\left(\displaystyle
1 - \alpha^{-\rho(t)} \right)^{-\gamma(t)/\rho(t)}
$$
and
$$
\Delta(\alpha,t) =
-\gamma(t)\alpha^{-\rho(t)} (1+O(\alpha^{-\rho(t)})),
$$
with $\rho(t)<0$. These results are collected in Table~\ref{tablex}.
In both Fr\'echet and Burr cases, $\Delta(\alpha,t)$ is asymptotically proportional to $\alpha^{-\rho(t)}$ as $\alpha \to 0$ with the convention $\rho(t)=-1$ for the Fr\'echet distribution.
Note that $\rho(t)$ is known as the second-order parameter in 
the extreme-value theory. It drives the quality of the approximation
of the conditional quantile $q(\alpha,t)$ by the power function $\alpha^{-\gamma(t)}$.
Furthermore, it is easily seen, that for these two distributions, the function $|\Delta(.,t)|$ is increasing. Thus, condition~(\ref{cond1}) of Corollary~\ref{coro}
can be simplified as 
$ m_t^{2\rho(t)} k_t^{1-2\rho(t)}\to 0$
which shows that, the smaller $\rho(t)$ is, the larger $k_t$ can be.
Finally, if $\gamma$ and $\rho$ are Lipschitzian, {\it i.e.} if there exist constants $c_{\gamma}>0$ and $c_{\rho}>0$ such that 
$$
|\gamma(x)-\gamma(x')| \leq c_{\gamma}d(x,x') \ {\rm{and}} \ |\rho(x)-\rho(x')| \leq c_{\rho}d(x,x')
$$
for all $(x,x')\in B(t,h_t)^2$,
then the oscillation can be bounded by $\omega_n(a)=O(h_t \log(1/a))$ as $a\to 0$ and thus
condition~(\ref{cond2}) of Corollary~\ref{coro} can be simplified as 
$ k_t^2 h_t \log m_t\to 0. $

\section{Finite sample behaviour}
\label{realdata}

\begin{figure}[p]
\begin{center}
\epsfig{figure=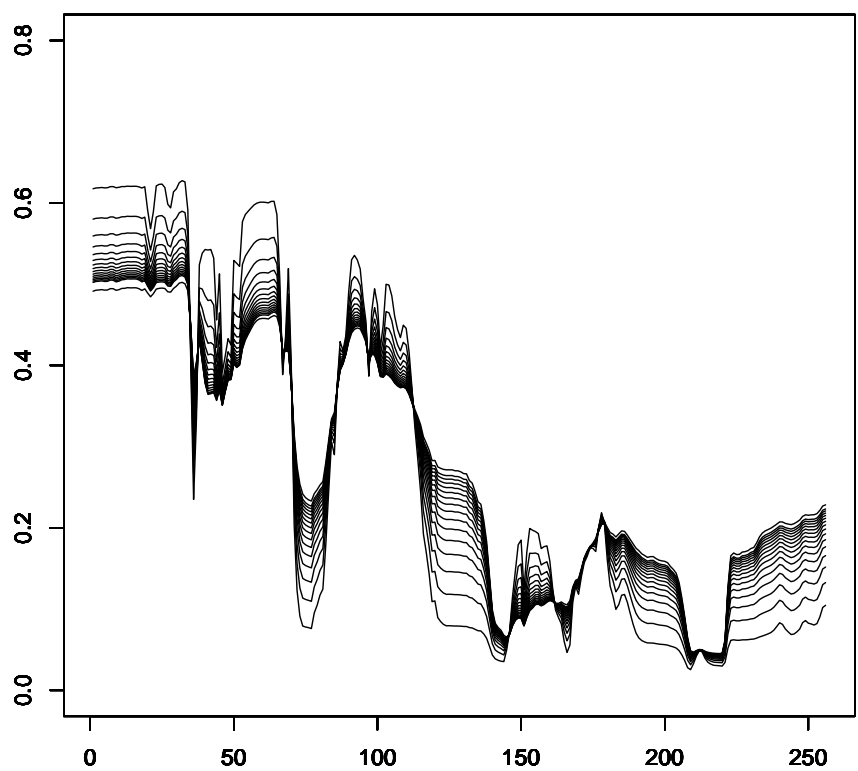,width=1\textwidth}
\caption{Representation of the 16 spectra as functions of the wavelength.}
\label{spectra}
\end{center}
\end{figure}

\begin{figure}[p]
\begin{center}
\epsfig{figure=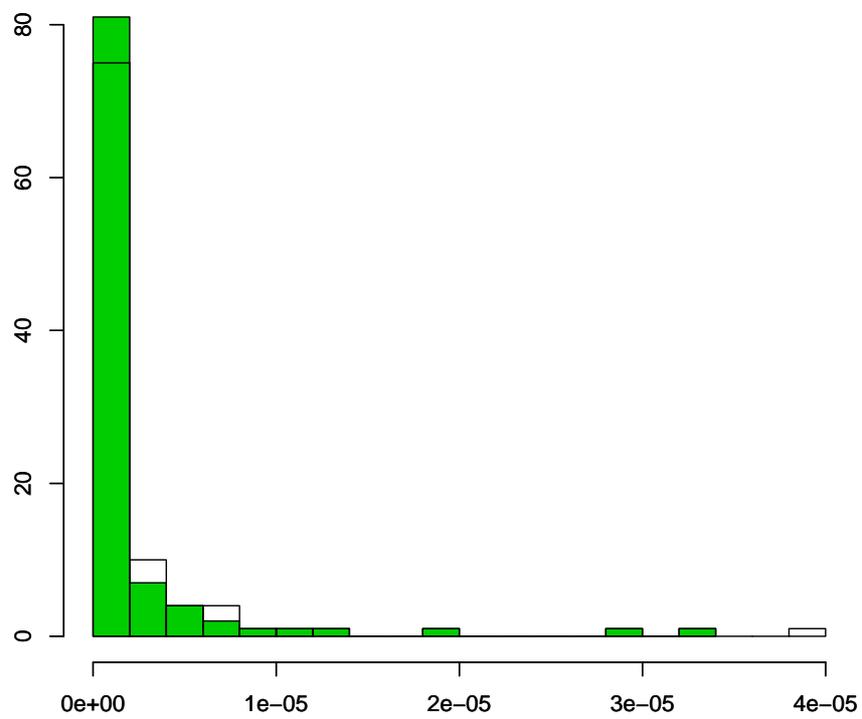,width=1\textwidth}
\caption{Comparison between the error distributions obtained with the heuristics method (transparent) and the oracle method (gray) on $N=100$ samples.}
\label{hist}
\end{center}
\end{figure}

\begin{figure}[p]
\begin{center}
\epsfig{figure=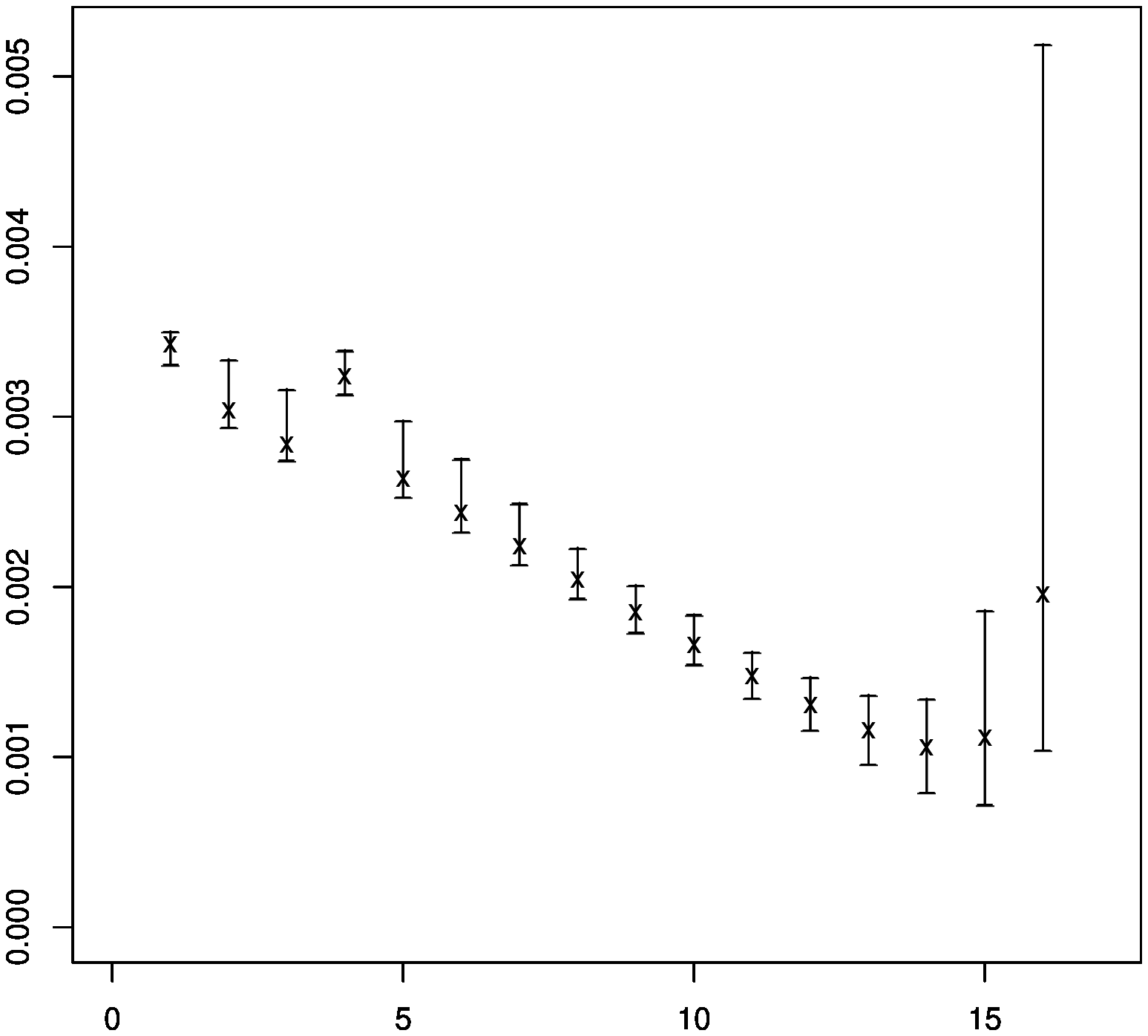,width=1\textwidth}
\caption{90\%-Empirical confidence intervals of $\hat q_2(1/300,.,\WZ)$ ranked by ascending order of the tail index.}
\label{boxplot1}
\end{center}
\end{figure}

\begin{figure}[p]
\begin{center}
\epsfig{figure=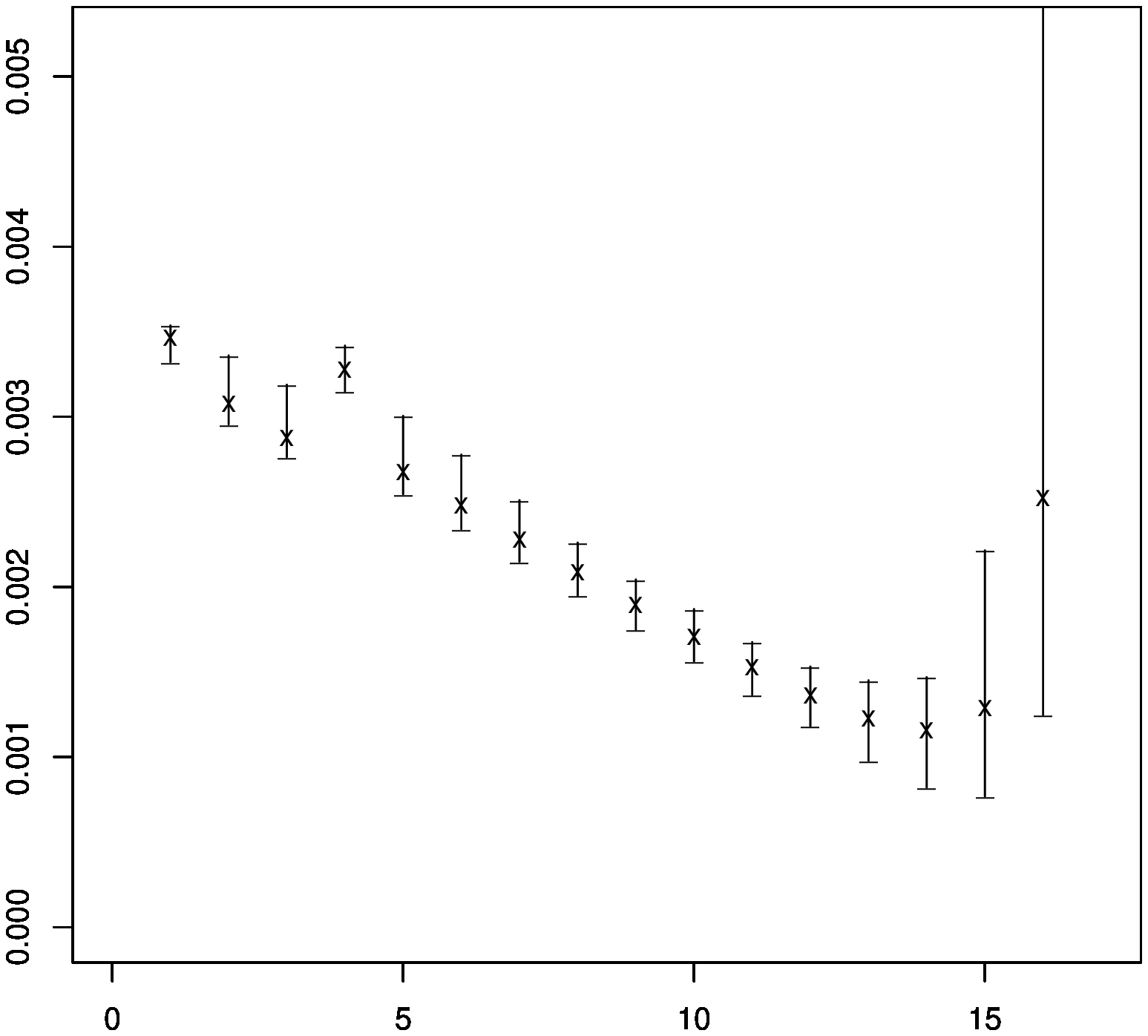,width=1\textwidth}
\caption{90\%-Empirical confidence intervals of $\hat q_2(1/500,.,\WZ)$ ranked by ascending order of the tail index.}
\label{boxplot2}
\end{center}
\end{figure}

\begin{figure}[p]
\begin{center}
\epsfig{figure=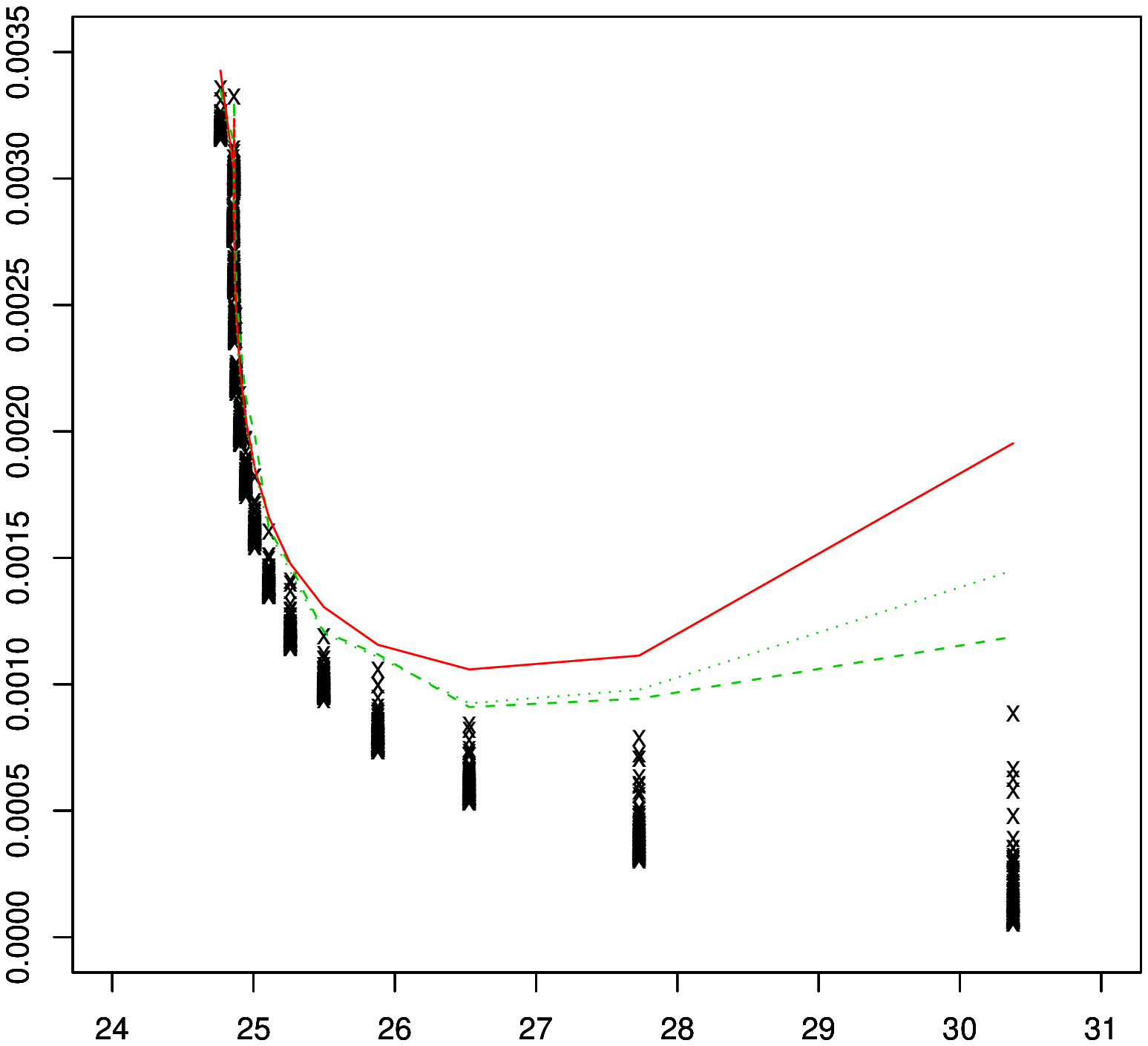,width=1\textwidth}
\caption{Comparison of the true quantile of order $\alpha=1/300$ (solid line) with the estimated quantiles by the heuristics strategy (dashed line) and the oracle strategy (dotted line) on the replication corresponding to the median error. The associated sample is represented by the points ("$\times$").}
\label{echmedian1}
\end{center}
\end{figure}

\begin{figure}[p]
\begin{center}
\epsfig{figure=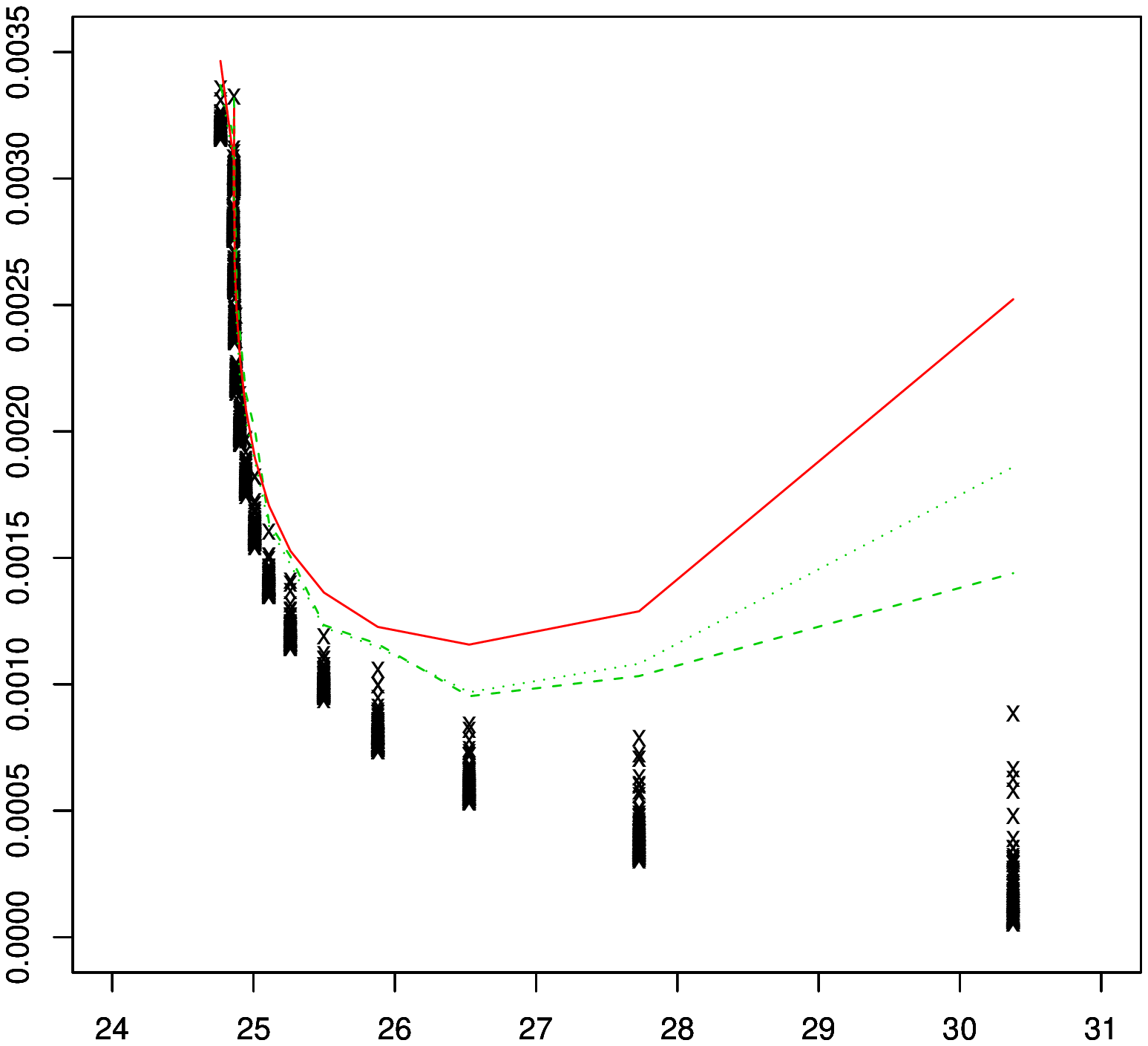,width=1\textwidth}
\caption{Comparison of the true quantile of order $\alpha=1/500$ (solid line) with the estimated quantiles by the heuristics strategy (dashed line) and the oracle strategy (dotted line) on the replication corresponding to the median error. The associated sample is represented by the points ("$\times$").}
\label{echmedian2}
\end{center}
\end{figure}

In this section, we propose to illustrate the behaviour of our conditional extreme quantiles estimators on functional spectrometric data. A question of interest for the planetologist is the following: Given a spectrum collected by the OMEGA instrument onboard the European spacecraft Mars Express in orbit around Mars, how to estimate the associated physical properties of the ground (grain size of CO$_2$, proportions of water, dust and CO$_2$, etc \ldots)?  To answer this question, a learning dataset can be constructed using radiative transfer models. Here, we focus on the CO$_2$ proportion. Given different values $y_i$, $i=1,\ldots,16$ of this proportion, a radiative transfer model provides us the corresponding spectra $x_i$, $i=1,\ldots,16$ (see Figure~\ref{spectra}). Clearly, the obtained spectra are non random. They are functions of the wavelength and we consider here their discretized version on 256 wavelengths $x_{i,l}$. Using this learning dataset, a lot of methods can be found in the literature to estimate the CO$_2$ proportion associated to an observed spectrum. One can mention Support Vector Machine, Sliced Inverse Regression, nearest neighbor approach, \ldots (see for instance~\cite{AGU} for an overview of these approaches). For all these methods, the estimation of the CO$_2$ proportion is perturbed by a random error term. We propose to modelize this perturbation by:
$$
Y_{i,j} = \log(1/y_i)+\sigma (\epsilon_j(x_i)-\Gamma(1-\gamma(x_i))), j=1,\ldots,n_i, \ i=1,\ldots,16,
$$
where
$$
\gamma(x_i)=0.3 \frac{\|x_i\|_2^2 - \min\limits_l \|x_l\|_2^2}{\max\limits_l \|x_l\|_2^2 - \min\limits_l \|x_l\|_2^2} + 0.2, \ \sigma=\min_i \frac{\log(1/y_i)}{\Gamma(1-\gamma(x_i))},
$$
and $\epsilon_j(x_i)$, $j=1,\ldots,n_i$ are independent and identically distributed random values from a Fr\'echet distribution with tail index $\gamma(x_i)$ (see Table~\ref{tablex}). Note that $\|x_i\|_2^2$ is an approximation of the total energy of the spectrum $x_i$. The above definitions ensure that $\gamma(x_i) \in [0.2,0.5]$ and that $Y_{i,j}>0$ for all $i=1,\ldots,16$ and $j=1,\ldots,n_i$. Furthermore, since the expectation of $\epsilon_j(x_i)$ is given by $\Gamma(1-\gamma(x_i))$, the random variables $Y_{i,j}$ are centered on the value $\log(1/y_i)$. Our aim is to estimate the conditional quantile
$$
q(\alpha,x_i)={\bar{F}}^{\leftarrow}(\alpha,x_i), \ {\rm{for}} \ i=1,\ldots,16,
$$
where ${\bar{F}}(.,x_i)$ is the survival distribution function of $Y_{i,1}$. To this end, the estimator
${\hat{q}}_2(\alpha,x_i,\WZ)$ defined in paragraph~\ref{exestis} is considered. The semi-metric distance based on the second derivative is 
adopted, as advised in~\cite{fer06}, Chapter~9:
$$
d^2(x_i,x_j)= \int \left ( x_i^{(2)}(t) - x_j^{(2)}(t) \right )^2dt,
$$
where $x^{(2)}$ denotes the second derivative of $x$. 
To compute this semi-metric, one can use an approximation of the functions $x_i$ and $x_j$ based on B-splines as proposed in~\cite{fer06}, Chapter~3. Here, we limit ourselves to a discretized version $\tilde d$ of $d$:
$$
\tilde d^2(x_i,x_j) = \sum_{l=2}^{255} \left\{ (x_{i,l+1} - x_{j,l+1}) + 
 (x_{i,l-1} - x_{j,l-1}) - 2 (x_{i,l}-x_{j,l}) \right\}^2.
$$
The finite sample performance of the estimator in assessed on $N=100$ replications of the sample $\{(x_i,Y_{i,j}), \ i=1,\ldots,16, \ j=1,\ldots,n_i\}$ with $n_1=\ldots=n_{16}=100$. Two values of $\alpha$ are considered: $1/300$ and $1/500$. In the following, we assume that the hyperparameters $h_t$ and $k_t$ does not depend on the spectrum (we thus omit the index $t$). These parameters are selected thanks to the heuristics proposed in~\cite{JMVA} which consists in minimizing the distance between two different estimators of the conditional extreme quantile:
$$
(\hat h_{{\rm{select}}},\hat k_{{\rm{select}}}) = \argmin_{h,k} \Delta({\hat{q}}_2(\alpha,.,\WH),{\hat{q}}_2(\alpha,.,\WZ),
$$
where for two functions $f$ and $g$,
$$
\Delta(f,g) = \left \{ \sum_{i=1}^{16} (f(x_i)-g(x_i))^2 \right \}^{1/2}.
$$
The estimator associated to these parameters is denoted by ${\hat{q}}_{{\rm{select}}}$. We also compute $\hat h_{\rm{oracle}}$ and $\hat k_{\rm{oracle}}$ defined as:
$$
(\hat h_{\rm{oracle}},\hat k_{\rm{oracle}}) = \argmin_{h,k} \Delta({\hat{q}}_2(\alpha,.,\WH),q(\alpha,.)).
$$
The conditional quantile estimator associated to these parameters is denoted by ${\hat{q}}_{{\rm{oracle}}}$. Note that $\hat h_{{\rm{select}}}$, $\hat k_{{\rm{select}}}$, $\hat h_{\rm{oracle}}$ and $\hat k_{\rm{oracle}}$ do not depend on $\alpha$. Of course, the oracle method cannot be applied in practical situations where $q(\alpha,.)$ is unknown. However, it provides us the lower bound on the distance $\Delta$ that can be reached with our estimator. In order to validate our choice of $\hat h_{{\rm{select}}}$ and $\hat k_{{\rm{select}}}$, the histograms of 
$
\Delta({\hat{q}}_{\rm{select}}(\alpha,.,\WZ),q(\alpha,.))$ and
$\Delta({\hat{q}}_{\rm{oracle}}(\alpha,.,\WZ),q(\alpha,.)),
$
computed for the $N=100$ replications, are superimposed on Figure~\ref{hist}. It appears that the mean errors are approximatively equal. Let us also remark that the heuristics errors seem to have a heavier right-tail than the oracle errors. For each spectrum $x_i$ the empirical 90\%-confidence interval of ${\hat{q}}_{{\rm{opt}}}(\alpha,x_i,\WZ)$ is represented on Figure~\ref{boxplot1} for $\alpha=1/300$ and on Figure~\ref{boxplot2} for $\alpha=1/500$. The confidence intervals are ranked by ascending order of the tail index. The larger the tail index is, the larger the confidence intervals are. This is in adequation with the result presented in Corollary~\ref{coro}. Finally, on Figure~\ref{echmedian1} ($\alpha=1/300$) and Figure~\ref{echmedian2} ($\alpha=1/500$), we draw estimators ${\hat{q}}_{{\rm{select}}}(\alpha,x_i,\WZ)$ and ${\hat{q}}_{{\rm{oracle}}}(\alpha,x_i,\WZ)$ as a function of $\|x_i\|_2^2$ on the replication giving rise to the median error $\Delta({\hat{q}}_{\rm{select}}(\alpha,.,\WZ),q(\alpha,.))$. It appears that the oracle estimator is only slightly better than the heuristics one. As noticed previously, the estimation error increases with the tail index.

\section{Proofs}
\label{secpreuves}

\subsection{Preliminary results}

Our first auxiliary lemma is a simple unconditioning tool for determining
the asymptotic distribution of a random variable.
\begin{Lem}
\label{decond}
Let $(X_n)$ and $(Y_n)$ be two sequences of real random variables. Suppose there exists an event ${\cal{A}}_n$ such that $(X_n|{\cal{A}}_n) \egloi (Y_n|{\cal{A}}_n)$ with $\proba({\cal{A}}_n) \to 1$. Then, $Y_n \cvloi Y$ implies $X_n \cvloi Y$.
\end{Lem}

\noindent {\bf{Proof of Lemma~\ref{decond}}} $-$ For all $x \in \R$,
the well-known expansion
\[ \proba(X_n \leq x) = \proba(\{X_n \leq x \}| {\cal{A}}_n)\proba({\cal{A}}_n) + \proba(\{X_n \leq x\} | {\cal{A}}_n^{C})\proba({\cal{A}}_n^{C}), \]
where ${\cal{A}}_n^{C}$ is the complementary event associated to ${\cal{A}}_n$, leads to the following inequalities:
\[ \proba(\{X_n \leq x\} | {\cal{A}}_n)\proba({\cal{A}}_n) \leq \proba(X_n \leq x) \leq \proba(\{X_n \leq x\} | {\cal{A}}_n)\proba({\cal{A}}_n) + \proba({\cal{A}}_n^{C}). \]
Since $(X_n|{\cal{A}}_n) \egloi (Y_n|{\cal{A}}_n)$, it follows that:
\[ \proba(\{Y_n \leq x \}\cap {\cal{A}}_n) \leq \proba(X_n \leq x) \leq \proba(\{Y_n \leq x\} \cap {\cal{A}}_n) + \proba({\cal{A}}_n^{C}). \]
Taking into account of
$$
\proba(Y_n \leq x) - \proba({\cal{A}}_n^{C})\leq\proba(\{Y_n \leq x\} \cap {\cal{A}}_n) \leq \proba(Y_n \leq x)
$$
leads to:
\[ \proba(Y_n \leq x) - \proba({\cal{A}}_n^{C}) \leq \proba(X_n \leq x) \leq \proba(Y_n \leq x) + \proba({\cal{A}}_n^{C}). \]
The conclusion is then straightforward since $\proba(Y_n \leq x) \to \proba(Y \leq x)$ and $\proba({\cal{A}}_n^{C}) \to 0$. \CQFD

\noindent The next lemma provides the asymptotic distribution of 
extreme quantile
estimators from an uniform distribution in a situation analogous to
{\bf (S.1)} in the unconditional case.
\begin{Lem}
\label{lemunif}
Let $V_1,\ldots,V_M$ be independent uniform random variables. For any sequence $(\theta_M) \subset (0,1)$ such that $\theta_M \to 0$ and $M\theta_M \to \infty$,
$$
\left ( \frac{M}{\theta_M} \right )^{1/2} (V_{\lfloor M\theta_M\rfloor,M}-\theta_M) \cvloi {\cal{N}}(0,1).
$$
\end{Lem}

\noindent {\bf{Proof of Lemma~\ref{lemunif}}} $-$ For the sake of simplicity,
let us introduce $k_M=\lfloor M\theta_M\rfloor$.
From R\'enyi's representation theorem,
$$
V_{k_M,M} \egloi \sum_{i=1}^{k_M} E_i \left / \sum_{i=1}^{M+1} E_i \right . 
$$
where $E_1,\ldots,E_{M+1}$ are independent random variables from a standard exponential distribution. Thus,
\begin{eqnarray*}
 \xi_{M} & \egdef & \left ( \frac{M}{\theta_M} \right )^{1/2} (V_{k_M,M}-\theta_M) \egloi \left ( \frac{1}{M} \sum_{i=1}^{M+1} E_i \right )^{-1} \left ( \frac{M}{\theta_M} \right )^{1/2} \\
 \ & \times & \left [ \frac{1}{k_M} \sum_{i=1}^{k_M} E_i \left ( \frac{k_M}{M} - \theta_M \right ) + \theta_M \left ( \frac{1}{k_M} \sum_{i=1}^{k_M} E_i-1\right ) \right . \\
 \ & - & \left .\theta_M \left ( \frac{1}{M} \sum_{i=1}^{M+1} E_i-1\right ) \right ],
\end{eqnarray*}
and, in view of the law of large numbers, we have 
\begin{eqnarray*}
\xi_M & \eqproba & \left ( \frac{M}{\theta_M} \right )^{1/2} \left ( \frac{k_M}{M} - \theta_M \right )(1+o_P(1)) +(M\theta_M)^{1/2} \left ( \frac{1}{k_M} \sum_{i=1}^{k_M} E_i-1\right ) \\
 \ & - & (M\theta_M)^{1/2} \left ( \frac{1}{M} \sum_{i=1}^{M+1} E_i-1\right ) \egdef \xi_{1,M} + \xi_{2,M} - \xi_{3,M}.
\end{eqnarray*}
Let us consider the three terms separately.
First, writing $k_M = M\theta_M - \tau_M$ with $\tau_M \in [0,1)$, we have
\begin{equation}
\label{T1N}
\xi_{1,M} \eqproba \left ( \frac{M}{\theta_M} \right )^{1/2} \frac{\tau_M}{M} = \frac{\tau_M}{(M\theta_M)^{1/2}} \to 0, 
\end{equation}
since $M\theta_M \to \infty$. Second, 
since $k_M \sim M\theta_M$, the central limit theorem entails
\begin{equation}
\label{T2N}
\xi_{2,M} \sim k_M^{1/2} \left ( \frac{1}{k_M} \sum_{i=1}^{k_M} E_i-1\right ) \cvloi {\cal{N}}(0,1).
\end{equation}
Similarly, it is easy to check that
\begin{equation}
\label{T3N}
\xi_{3,M} = O_{\rm{P}} (\theta_M^{1/2}) = o_{\rm{P}}(1), 
\end{equation}
since $\theta_M \to 0$. Collecting (\ref{T1N}), (\ref{T2N}) and (\ref{T3N}) concludes the proof. \CQFD

\subsection{Proofs of main results}
\label{mainresults}

{\bf{Proof of Proposition~\ref{proplaurent}}} $-$ 
Under {\bf (A)} and
since the random values $\{ Z_{i}(t), \ i=1,\ldots,m_t\}$ are independent, we have:
\[ \{ \log Z_i(t), \ i=1,\ldots,m_t\} \egloi \{ \log q(V_i,x_i)  \ i=1,\ldots,m_t \}, \]
where $x_i$ is the covariate associated to $Z_i(t)$. Denoting by $\psi(i)$ the random index of the covariate associated to the observation $Z_{m_t-i+1,m_t}(t)$,
we obtain
\[ \{ \log Z_{m_t-i+1,m_t}(t), \ i=1,\ldots,m_t\} \egloi \{ \log q(V_{\psi(i)},x_{\psi(i)})  \ i=1,\ldots,m_t \}. \]
Let us consider the event $ {\cal{A}}_n = {\cal{A}}_{1,n} \cap {\cal{A}}_{2,n}$
where 
\begin{eqnarray*}
 {\cal{A}}_{1,n} &=& \left \{ \min_{i=1,\ldots,k_t-1} \log \frac{q(V_{i,m_t},u_{i})}{q(V_{i+1,m_t},u_{i+1})} > 0, \forall (u_{1},\ldots,u_{k_t}) \subset B(t,h_t) \right \} 
\mbox{ and }\\
 {\cal{A}}_{2,n} &=& 
\left \{ \min_{i=k_t+1,\ldots,m_t} \log \frac{q(V_{k_t,m_t},u_{k_t})}{q(V_{i,m_t},u_{i})} > 0, \forall (u_{k_t+1},\ldots,u_{m_t}) \subset B(t,h_t) \right \}.
\end{eqnarray*}
Conditionally to $ {\cal{A}}_{1,n} $, the random variables 
$q(V_{i,m_t},u_{i})$, $i=1,\dots, k_t$ are ordered as
$$
q(V_{k_t,m_t},u_{k_t})\leq q(V_{k_t-1,m_t},u_{k_t-1})
\leq \dots \leq q(V_{1,m_t},u_{1}),
$$
and, conditionally to $ {\cal{A}}_{2,n} $, the remaining random variables
$q(V_{i,m_t},u_{i})$, $i=k_t+1,\dots, m_t$ 
are smaller since
$$
\max_{i=k_t+1,\ldots,m_t}q(V_{i,m_t},u_{i}) \leq q(V_{k_t,m_t},u_{k_t}).
$$
Thus, conditionally to ${\cal{A}}_n$, the $k_t$ largest random values taken from the set $\{ \log q(V_{\psi(i)},x_{\psi(i)}), \ i=1,\ldots,m_t \}$ are $\{ \log q(V_{i,m_t},x_{\psi(i)}), \ i=1,\ldots,k_t \}$. 
Consequently, for $J_{k_t}=\{ 1,\ldots,k_t \}$ and letting $T_i\egdef x_{\psi(i)}$, we have:
\[ \left \{ \log Z_{m_{t}-i+1,m_{t}}(t), i \in J_{k_t} |  {\cal{A}}_n \right \} \egloi \left \{ \log q(V_{i,m_{t}},T_i), i \in J_{k_t} |  {\cal{A}}_n \right \}. \]
To conclude the proof, it remains to show that $\proba({\cal{A}}_n) \to 1$ as $n\to\infty$. 
Let us define $\delta_{m_t}=m_t^{-(1+\delta)}$ and consider
the events 
\begin{eqnarray*}
{\cal{A}}_{3,n} &=& \{ V_{1,m_t} > \delta_{m_t}\} \cap \{V_{m_t,m_t} < 1-\delta_{m_t}\}\\
{\cal{A}}_{4,n} &=& \left \{ \min_{i=1,\ldots,k_t} \log \frac{q(V_{i,m_t},t)}{q(V_{i+1,m_t},t)} > 2\omega_n(\delta_{m_t}) \right \}. 
\end{eqnarray*}
Under ${\cal{A}}_{3,n}$, we have
$\delta_{m_t} < V_{i,m_t} < 1-\delta_{m_t}$
for all $i=1,\ldots,m_t$.
Hence, for all $(u_i,u_j) \in B(t,h_t)^2$, it follows that, on the one hand 
\begin{eqnarray*} 
\log \frac{q(V_{j,m_t},u_j)}{q(V_{i,m_t},u_i)} & = &  \log \frac{q(V_{j,m_t},t)}{q(V_{i,m_t},t)} + \log \frac{q(V_{j,m_t},u_j)}{q(V_{j,m_t},t)} + \log \frac{q(V_{i,m_t},t)}{q(V_{i,m_t},u_i)} \\
 \ & \geq & \log \frac{q(V_{j,m_t},t)}{q(V_{i,m_t},t)} - 2\omega_n(\delta_{m_t}),
\end{eqnarray*}
and on the other hand,
\begin{eqnarray*}
\min_{i=k_t+1,\ldots,m_t} \log \frac{q(V_{k_t,m_t},u_{k_t})}{q(V_{i,m_t},u_i)} & \geq & \min_{i=k_t+1,\ldots,m_t} \log \frac{q(V_{k_t,m_t},t)}{q(V_{i,m_t},t)} - 2\omega_n(\delta_{m_t}) \\
 \ & \geq & \log \frac{q(V_{k_t,m_t},t)}{q(V_{k_t+1,m_t},t)} - 2\omega_n(\delta_{m_t}).
\end{eqnarray*}
Consequently ${\cal{A}}_{3,n} \cap {\cal{A}}_{4,n} \subset {\cal{A}}_n$. 
Remarking that
$$
\proba({\cal{A}}_{3,n}) \geq \proba(V_{1,m_t} > \delta_{m_t})+\proba(V_{m_t,m_t} < 1-\delta_{m_t})-1 = 2\proba(V_{1,m_t} > \delta_{m_t})-1 \to 1,
$$
since $V_{m_t,m_t} \egloi 1-V_{1,m_t}$ and $\proba ( V_{1,m_t}>\delta_{m_t}) = \left ( 1-\delta_{m_t} \right )^{m_t} \to 1$, it thus remains to prove that $\proba({\cal{A}}_{4,n}) \to 1$. 
From \cite{BING}, paragraph~1.3.1, condition {\bf{(A)}} implies that there exists $c(t)>0$, depending only on $t$ such that, for all $\alpha \in (0,1)$,
$$
q(\alpha,t) = c(t) \exp \left \{ \int_{\alpha}^1 \frac{\gamma(t)+\Delta(u,t)}{u}du\right \},
$$
which is the so-called Karamata representation for normalised regularly varying functions.
Hence, for all $i \in J_{k_t}$, 
$$
\log \frac{q(V_{i,m_t},t)}{q(V_{i+1,m_t},t)} = 
\int_{V_{i,m_t}}^{V_{i+1,m_t}} \frac{\gamma(t)+\Delta(u,t)}{u} du,
$$
and it follows that
$$
\log \frac{q(V_{i,m_t},t)}{q(V_{i+1,m_t},t)} \geq (\gamma(t)- {\bar{\Delta}}(V_{k_t+1,m_t},t)) \log \frac{V_{i+1,m_t}}{V_{i,m_t}},
$$
leading to
\begin{eqnarray*}
\label{decompo}
\proba({\cal{A}}_{4,n}) & \geq  & \proba \left( (\gamma(t)- {\bar{\Delta}}(V_{k_t+1,m_t},t))  
 \min_{i=1,\ldots,k_t}  \log \frac{V_{i+1,m_t}}{V_{i,m_t}} > 2 \omega_n(\delta_{m_t}) \right ) \\
 \ & \geq & \proba \left ( \left \{ \min_{i=1,\ldots,k_t} 
 \log \frac{V_{i+1,m_t}}{V_{i,m_t}} \geq \frac{4\omega_n(\delta_{m_t})}{\gamma(t)} \right \} \cap \left \{ {\bar{\Delta}}(V_{k_t+1,m_t},t) < \gamma(t)/2 \right \} \right )  \\
 \ & \geq & \proba \left (  \min_{i=1,\ldots,k_t} 
 \log \frac{V_{i+1,m_t}}{V_{i,m_t}} \geq \frac{4\omega_n(\delta_{m_t})}{\gamma(t)} \right) + \proba\left( {\bar{\Delta}}(V_{k_t+1,m_t},t) < \gamma(t)/2 \right) -1 \\
 \ & \egdef & P_{1,m_t} + P_{2,m_t}-1.
\end{eqnarray*}
In view of R\'enyi representation for uniform ordered random variables,
$$
\{i\log(V_{i,m_t}^{-1}/V_{i+1,m_t}^{-1}), \ i\in J_{k_t} \} \egloi \{ F_i, \ i\in J_{k_t} \} ,
$$
where $F_1,\ldots,F_{k_t}$ are independent random variables from a standard exponential distribution, we have
\begin{eqnarray*}
P_{1,m_t} &=& \proba \left ( \min_{i=1,\ldots,k_t} \frac{F_i}{i} \geq \frac{4\omega_n(\delta_{m_t})}{\gamma(t)}  \right ) 
= \prod_{i=1}^{k_t} \exp\left(- \frac{4i\omega_n(\delta_{m_t})}{\gamma(t)} \right)\\
 & = & \exp \left ( -  \frac{2}{\gamma(t)} k_t(k_t+1) \omega_n(\delta_{m_t}) \right ) 
 \to 1, 
\end{eqnarray*}
since $k_t^2\omega_n(\delta_{m_t})\to 0$. Furthermore, 
$V_{k_t+1,m_t} = (k_t/m_t)(1+o_{\rm{P}}(1)) \cvproba 0$ and $\Delta(\alpha,t) \to 0$ as $\alpha\to 0$
entail $P_{2,m_t}\to 1$. The conclusion follows.\CQFD

\noindent {\bf{Proof of Proposition~\ref{compamax}}} $-$ From Proposition \ref{proplaurent}, there exists an event ${\cal{A}}_n$ with $\proba({\cal{A}}_n) \to 1$ such that $(Z_{m_t,m_t}(t)|{\cal{A}}_n) \egloi  (q(V_{1,m_t},T_1)|{\cal{A}}_n)$ and thus,
\begin{eqnarray} 
\label{decompo2}
\proba(Z_{m_t,m_t}(t)<q(\alpha_{m_t},t)) & = & \proba \left ( \left\{\log \frac{q(V_{1,m_t},T_1)}{q(\alpha_{m_t},t)}<0 \right\}\cap {\cal{A}}_n \right ) \nonumber \\
 \ & + & \proba \left (\left\{\log \frac{Z_{m_t,m_t}(t)}{q(\alpha_{m_t},t)}<0 \right\}\cap {\cal{A}}_n^{C} \right ) \nonumber \\
 \ & \egdef & P_{3,m_t}+P_{4,m_t}.
\end{eqnarray}
Clearly, $P_{4,m_t} \leq \proba({\cal{A}}_n^{C}) \to 0$. Let us now consider the term $P_{3,m_t}$. Introducing $\delta_{m_t}=m_t^{-(1+\delta)}$ and
${\cal{A}}_{5,n}=\{ V_{1,m_t} \in [\delta_{m_t},1-\delta_{m_t}] \}$, we have
\begin{eqnarray*} 
P_{3,m_t} & = & \proba \left ( \left\{\log \frac{q(V_{1,m_t},T_1)}{q(\alpha_{m_t},t)}<0 \right\}\cap {\cal{A}}_n \cap {\cal{A}}_{5,n}  \right ) \\
 \ & + & \proba \left (\left\{\log \frac{q(V_{1,m_t},T_1)}{q(\alpha_{m_t},t)}<0 \right\}\cap {\cal{A}}_n \cap {\cal{A}}_{5,n}^C  \right )
\end{eqnarray*}
and standard calculations lead to:
\begin{eqnarray*}
 &  & \proba \left (\left\{\log \frac{q(V_{1,m_t},T_1)}{q(\alpha_{m_t},t)}<0\right\} \cap {\cal{A}}_{5,n} \right ) + \proba({\cal{A}}_n) - 1 \leq  P_{3,m_t}\\ 
  \leq &&\proba \left (\left\{\log \frac{q(V_{1,m_t},T_1)}{q(\alpha_{m_t},t)}<0 \right\} \cap {\cal{A}}_{5,n} \right ) + \proba({\cal{A}}_{5,n}^C ).
\end{eqnarray*}
Furthermore, ${\cal{A}}_{5,n}$ implies
\[ \left| \log \frac{q(V_{1,m_t},T_1)}{q(V_{1,m_t},t)} \right|\leq \omega_n(\delta_{m_t}), \]
and thus
\begin{eqnarray*}
 &  & \proba \left (\left\{\log \frac{q(V_{1,m_t},t)}{q(\alpha_{m_t},t)} < - \omega_n(\delta_{m_t})\right\} \cap {\cal{A}}_{5,n} \right ) + \proba({\cal{A}}_n) - 1  \leq  P_{3,m_t}\\
 \leq&& \proba \left  (\left\{\log \frac{q(V_{1,m_t},t)}{q(\alpha_{m_t},t)} < \omega_n(\delta_{m_t}) \right\}\cap {\cal{A}}_{5,n} \right ) + \proba({\cal{A}}_{5,n}^C),
\end{eqnarray*}
which entails
\begin{eqnarray}
\label{cadre}
 &  & \proba\left ( \log \frac{q(V_{1,m_t},t)}{q(\alpha_{m_t},t)} < - \omega_n(\delta_{m_t}) \right ) + \proba ( {\cal{A}}_{5,n}) + \proba({\cal{A}}_n) - 2 \leq  P_{3,m_t} \nonumber \\
  \leq &&\proba \left ( \log \frac{q(V_{1,m_t},t)}{q(\alpha_{m_t},t)} < \omega_n(\delta_{m_t}) \right ) + \proba({\cal{A}}_{5,n}^C).
\end{eqnarray}
Let us now focus on the quantity
\begin{eqnarray*}
 P_{5,m_t} &\egdef& \proba\left ( \log \frac{q(V_{1,m_t},t)}{q(\alpha_{m_t},t)} < \pm \omega_n(\delta_{m_t}) \right )\\
& = &\left [ \proba \left ( \log \frac{q(V_{1},t)}{q(\alpha_{m_t},t)} < \pm \omega_n(\delta_{m_t}) \right ) \right ]^{m_t}\\
& =& \left [ \proba \left ( q(V_{1},t) < {\rm{e}}^{\pm \omega_n(\delta_{m_t})} q(\alpha_{m_t},t) \right ) \right ]^{m_t}\\
 & = & \left [ \proba \left ( 1-V_{1} < F \left( {\rm{e}}^{\pm \omega_n(\delta_{m_t})} q(\alpha_{m_t},t),t \right ) \right ) \right ]^{m_t}\\ 
 \ & = & \exp \left [ m_t \log F \left ( {\rm{e}}^{\pm \omega_n(\delta_{m_t})} q(\alpha_{m_t},t),t \right ) \right ].
\end{eqnarray*}
Since ${\rm{e}}^{\pm \omega_n(\delta_{m_t})} q(\alpha_{m_t},t) \to \infty$ and introducing the conditional survival function ${\bar{F}}(.,t)=1-F(.,t)$, we have
\begin{eqnarray*}
m_t \log F \left ( {\rm{e}}^{\pm \omega_n(\delta_{m_t})} q(\alpha_{m_t},t),t \right ) & = & -m_t {\bar{F}} \left ( {\rm{e}}^{\pm \omega_n(\delta_{m_t})} q(\alpha_{m_t},t),t \right ) (1+o(1)) \\
 \ & = & -m_t \alpha_{m_t} \frac{{\bar{F}} \left ( {\rm{e}}^{\pm \omega_n(\delta_{m_t})} q(\alpha_{m_t},t),t \right )}{ {\bar{F}} \left ( q(\alpha_{m_t},t),t \right )}(1+o(1)).
\end{eqnarray*}
As already mentioned, {\bf (A)} implies~(\ref{modelq}) 
which, in turn, shows that ${\bar{F}}(.,t)$ is a regularly function at infinity with index $-1/\gamma(t)$. Hence, since ${\rm{e}}^{\pm \omega_n(\delta_{m_t})} \to 1$, we thus have (see \cite{BING}, Theorem~1.5.2),
\[ \frac{{\bar{F}} \left ( {\rm{e}}^{\pm \omega_n(\delta_{m_t})} q(\alpha_{m_t},t),t \right )}{ {\bar{F}} \left ( q(\alpha_{m_t},t),t \right )} \to 1. \]
As a conclusion,
\begin{equation}
\label{alpha}
P_{5,m_t} = \left [ 1-\alpha_{m_t}(1+o(1)) \right ]^{m_t},
\end{equation}
and collecting (\ref{cadre}) and (\ref{alpha}) leads to:
\begin{eqnarray*}
 \ & \ & \left [ 1-\alpha_{m_t}(1+o(1)) \right ]^{m_t} + \proba ( {\cal{A}}_{5,n}) + \proba({\cal{A}}_n) - 2 \nonumber \\
 \ & \leq & P_{3,m_t} \leq \left [ 1-\alpha_{m_t}(1+o(1)) \right ]^{m_t} + \proba({\cal{A}}_{5,n}^C).
\end{eqnarray*}
Since $\proba({\cal{A}}_{5,n}) \to 1$ and $\proba({\cal{A}}_n) \to 1$, it is then straightforward that $P_{3,m_t} \to 0$ under {\bf{(S.1)}} and  $P_{3,m_t} \to {\rm{e}}^{-c}$ under {\bf{(S.2)}} or {\bf{(S.3)}}. Equation~(\ref{decompo2}) concludes the proof. \CQFD

\noindent {\bf{Proof of Theorem \ref{quantclass}}} $-$ 
Let us introduce, for the sake of simplicity,  $k_t=\lfloor m_t \alpha_{m_t}\rfloor$.
From Proposition~\ref{proplaurent}, there exists an event ${\cal{A}}_n$ such that:
\[ \left ( \left . (m_t\alpha_{m_t})^{1/2} \log \frac{{\hat{q}}_1(\alpha_{m_t},t)}{q(\alpha_{m_t},t)} \right | {\cal{A}}_n \right ) \egloi \left ( \left . (m_t\alpha_{m_t})^{1/2} \log \frac{q(V_{k_t,m_t},T_{k_t})}{q(\alpha_{m_t},t)} \right | {\cal{A}}_n \right ), \]
where $\proba({\cal{A}}_n) \to 1$. From~Lemma \ref{decond}, 
the convergence in distribution
\begin{equation}
\label{cvloiunif}
(m_t\alpha_{m_t})^{1/2} \log \frac{q(V_{k_t,m_t},T_{k_t})}{q(\alpha_{m_t},t)} \cvloi {\cal{N}}(0,\gamma^2(t)), 
\end{equation}
is a sufficient condition to obtain
\[ (m_t\alpha_{m_t})^{1/2} \log \frac{{\hat{q}}_1(\alpha_{m_t},t)}{q(\alpha_{m_t},t)} \cvloi {\cal{N}}(0,\gamma^2(t)). \]
A straightforward application of the $\delta$-method will then conclude the proof. Let us prove the convergence in distribution~(\ref{cvloiunif}). 
To this end, consider
\[ R_n = \left|\log \frac{q(V_{k_t ,m_t},T_{k_t})}{q(V_{k_t ,m_t},t)}\right| \]
and let $\delta_{m_t}=m_t^{-(1+\delta)}$. Remark that, under {\bf{(S.1)}},
$$
\proba(R_n \leq \omega_n(\delta_{m_t}))  \geq 
    \proba(V_{k_t ,m_t} \in [\delta_{m_t},1- \delta_{m_t}]) \to 1.
$$
Thus, $R_n=O_{\rm{P}}(\omega_n(\delta_{m_t}))$
and we have
\begin{equation}
\label{firstpart}
\log \frac{q(V_{k_t,m_t},T_{k_t} )}{q(\alpha_{m_t},t)} = \log \frac{q(V_{k_t,m_t},t)}{q(\alpha_{m_t},t)} + O_{\rm{P}}(\omega_n(\delta_{m_t})).
\end{equation}
Let us introduce the log-quantile function $g(.)=\log q(.,t)$. Clearly, for all
$\alpha \in (0,1)$,
$$
g'(\alpha) = \frac{\Delta(\alpha,t)-\gamma(t)}{\alpha}
$$
and a first-order Taylor expansion leads to:
\begin{eqnarray*}
(m_t\alpha_{m_t})^{1/2} \log \frac{q(V_{k_t ,m_t},t)}{q(\alpha_{m_t},t)} & = & (m_t\alpha_{m_t})^{1/2} g'(\theta_{m_t}) (V_{k_t,m_t}-\alpha_{m_t}) \\
 \ & = & \alpha_{m_t} g'(\theta_{m_t}) \left ( \frac{m_t}{\alpha_{m_t}} \right )^{1/2} (V_{ k_t ,m_t}-\alpha_{m_t}),
\end{eqnarray*}
where $\theta_{m_t} \in [\min(\alpha_{m_t},V_{k_t ,m_t}),\max(\alpha_{m_t},V_{k_t ,m_t})]$. Now, $V_{k_t,m_t} \eqproba \alpha_{m_t}$ entails $\theta_{m_t} \eqproba \alpha_{m_t} \to 0$ and, from {\bf{(A)}},
$$
\alpha_{m_t}g'(\theta_{m_t}) \eqproba \theta_{m_t}g'(\theta_{m_t}) = \Delta(\theta_{m_t},t) - \gamma(t) \cvproba -\gamma(t). 
$$
Then, Lemma~\ref{lemunif} implies that
\begin{equation}
\label{secondpart}
(m_t\alpha_{m_t})^{1/2} \log \frac{q(V_{k_t,m_t},t)}{q(\alpha_{m_t},t)} \cvloi {\cal{N}}(0,\gamma(t)^2). 
\end{equation}
Collecting (\ref{firstpart}) and (\ref{secondpart}) concludes the proof after remarking that condition $(m_t\alpha_{m_t})^{2}\omega_n(\delta_{m_t}) \to 0$ implies $(m_t\alpha_{m_t})^{1/2}\omega_n(\delta_{m_t}) \to 0$. \CQFD

\noindent {\bf{Proof of Theorem \ref{quantinter}}} $-$ Since $q(.,t)$ is regularly varying with index $-\gamma(t)$,
we have under {\bf{(S.2)}} that
$q(1/m_t,t)/q(\alpha_{m_t},t) \sim (m_t\alpha_{m_t})^{\gamma(t)} \to c^{\gamma(t)}$
and the following asymptotic expansion holds
\begin{eqnarray*}
\log \frac{{\hat{q}}_1(\alpha_{m_t},t)}{q(\alpha_{m_t},t)}
&=& \log \frac{{\hat{q}}_1(\alpha_{m_t},t)}{q(1/m_t,t)}
+ \log \frac{q(1/m_t,t)}{q(\alpha_{m_t},t)}\\
&=& \log \frac{{\hat{q}}_1(\alpha_{m_t},t)}{q(1/m_t,t)} + \gamma(t) \log(c) +o(1).
\end{eqnarray*}
Now, recall that in situation {\bf (S.2)}, for $n$ large enough,
$\lfloor m_t \alpha_{m_t}\rfloor = \lfloor c \rfloor$. Thus, from Proposition~\ref{proplaurent}, there exists an event ${\cal{A}}_n$ such that $\proba({\cal{A}}_n) \to 1$ and
\[ \left ( \left . \log \frac{{\hat{q}}_1(\alpha_{m_t},t)}{q(1/m_t,t)} \right | {\cal{A}}_n \right ) \egloi \left ( \left . \log \frac{q(V_{\lfloor c\rfloor,m_t},T_{\lfloor c\rfloor})}{q(1/m_t,t)} \right | {\cal{A}}_n \right ). \]
Mimicking the proof of Theorem~\ref{quantclass}, we obtain 
$$
\log \frac{q(V_{\lfloor c\rfloor,m_t},T_{\lfloor c\rfloor})}{q(1/m_t,t)} = \log \frac{q(V_{\lfloor c\rfloor,m_t},t)}{q(1/m_t,t)} + O_{\rm{P}}(\omega_n(\delta_{m_t})).
$$
To conclude, one can remark that $q(V_{\lfloor c\rfloor,m_t},t)$ is the
$\lfloor c\rfloor$th uppest order statistics associated 
to a heavy-tailed distribution. In such a case, Corollary~4.2.4 of~\cite{EMBR}
states that $q(V_{\lfloor c\rfloor,m_t},t)/q(1/m_t,t)$ converges to
a non-degenerated distribution. This asymptotic distribution is explicit
even though it is not reproduced here. \CQFD

\noindent {\bf{Proof of Theorem \ref{quantextreme}}} $-$ Observing that
$$
\log {\hat{q}}_2(\alpha_{m_t},t) = \log {\hat{q}}_1(\beta_{m_t},t) + {\hat{\gamma}}_n(t) \log \left ( \frac{\beta_{m_t}}{\alpha_{m_t}} \right )
$$
leads to the following expansion
\begin{eqnarray*}
\log \frac{{\hat{q}}_2(\alpha_{m_t},t)}{q(\alpha_{m_t},t)} & = & \log \frac{{\hat{q}}_1(\beta_{m_t},t)}{q(\beta_{m_t},t)} \\
 \ & + & \log \left ( \frac{\beta_{m_t}}{\alpha_{m_t}} \right ) ({\hat{\gamma}}_n(t)-\gamma(t) ) \\
 \ & - & \log \frac{q(\alpha_{m_t},t)}{q(\beta_{m_t},t)} + \gamma(t) \log \left ( \frac{\beta_{m_t}}{\alpha_{m_t}} \right ) \\
 \ & \egdef & \xi_{4,m_t} + \xi_{5,m_t} - \xi_{6,m_t}.
\end{eqnarray*}
First remark that, under {\bf{(A)}}, as already shown in the proof of Proposition~\ref{proplaurent},
$$
\log \frac{q(\alpha_{m_t},t)}{q(\beta_{m_t},t)} = \int_{\alpha_{m_t}}^{\beta_{m_t}} \frac{\gamma(t)+\Delta(u,t)}{u} du, 
$$
and thus, $\xi_{6,m_t}$ can be simplified as
$$
\xi_{6,m_t} = \int_{\alpha_{m_t}}^{\beta_{m_t}} \frac{\Delta(u,t)}{u} du
$$
which leads to the bound:
$$
|\xi_{6,m_t}|\leq {\bar{\Delta}}(\beta_{m_t},t)  \log \left ( \frac{\beta_{m_t}}{\alpha_{m_t}} \right ).
$$
The two additional conditions are now treated separately since, under condition~(\ref{C1}), the asymptotic distribution is imposed by $\xi_{4,m_t}$ whereas, under~(\ref{C2}), the asymptotic distribution is imposed by $\xi_{5,m_t}$.
 \ \\
\noindent (i) Under (\ref{C1}), Theorem~\ref{quantclass} entails that
\begin{equation}
\label{I1}
(m_t\beta_{m_t})^{1/2}\xi_{4,m_t} \cvloi {\cal{N}}(0,\gamma^2(t))
\end{equation}
and
\begin{equation}
\label{I2}
(m_t\beta_{m_t})^{1/2}\xi_{5,m_t} = \zeta_{m_t} \upsilon_n^{-1}(t)  \upsilon_n(t) ( {\hat{\gamma}}_n(t) - \gamma(t)) \cvproba 0,
\end{equation}
from (\ref{CVgamma}) and (\ref{C1}). Finally,
\begin{equation}
\label{I3}
(m_t\beta_{m_t})^{1/2}|\xi_{6,m_t}| \leq \zeta_{m_t} {\bar{\Delta}}(\beta_{m_t},t)  \to 0,
\end{equation}
from (\ref{C1}). Collecting (\ref{I1}), (\ref{I2}) and (\ref{I3}) concludes the proof of (\ref{cvloi1}). \\
 \ \\
\noindent (ii) Under (\ref{C2}),  Theorem~\ref{proplaurent} implies
\begin{equation}
\label{J1}
\frac{\upsilon_n(t)}{\log(\beta_{m_t}/\alpha_{m_t})} \xi_{4,m_t} = \upsilon_n(t)
\zeta_{m_t}^{-1}(m_t\beta_{m_t})^{1/2}\xi_{4,m_t} \cvproba 0.
\end{equation}
Moreover, from (\ref{CVgamma}), 
\begin{equation}
\label{J2}
\frac{\upsilon_n(t)}{\log(\beta_{m_t}/\alpha_{m_t})}\xi_{5,m_t} =\upsilon_n(t) ( {\hat{\gamma}}_n(t) - \gamma(t)) \cvloi {\cal{D}}
\end{equation}
and finally,
\begin{equation}
\label{J3}
\frac{\upsilon_n(t)}{\log(\beta_{m_t}/\alpha_{m_t})} |\xi_{6,m_t}| \leq {\bar{\Delta}}(\beta_{m_t},t) \upsilon_n(t) \to 0,
\end{equation}
under (\ref{C2}). Collecting (\ref{J1}), (\ref{J2}) and (\ref{J3}) concludes the proof of (\ref{cvloi2}). \CQFD

\begin{table}[p]
\begin{center}
$
\begin{array}{|l|l|l|}
\hline
&&\\
			& q(\alpha,t)  & \Delta(\alpha,t) \\
&&\\
\hline
&&\\
\mbox{Pareto} 	& {\alpha}^{-\gamma(t)} & 0 \\
&&\\
\mbox{Fr\'echet}	& {\alpha}^{-\gamma(t)} \left\{\displaystyle\frac{1}{\alpha}\log\left(\displaystyle
\frac{1}{1-\alpha}\right)\right\}^{-\gamma(t)}  & -\displaystyle\frac{\gamma(t)}{2} \alpha (1+O(\alpha)) \\
&&\\
\mbox{Burr}	& \alpha^{-\gamma(t)}\left(\displaystyle
1 - \alpha^{-\rho(t)} \right)^{-\gamma(t)/\rho(t)}&   -\gamma(t)\alpha^{-\rho(t)} (1+O(\alpha^{-\rho(t)}))   \\
&&\\
\hline
\end{array}
$
\end{center}
\caption{Some examples of heavy-tailed distributions. For all distributions, $\gamma(t)>0$ is the tail-index and $\rho(t)<0$ is referred to as the second-order parameter in extreme-value theory.}
\label{tablex}
\end{table}

\end{document}